\documentclass[twoside]{amsart}

\usepackage{latexsym, amsmath, amssymb, amsfonts, amscd, amsthm}

\newcommand{\Ce}{{\mathbb C}}

\renewcommand{\Re}{{\mathbb R}}

\newcommand{\Ze}{{\mathbb Z}}

\newcommand{\pa}{\parallel}

\newtheorem{theorem}{Theorem}[section]
\newtheorem{lemma}[theorem]{Lemma}
\newtheorem{proposition}[theorem]{Proposition}

\newtheorem{example}[theorem]{Example}
\newtheorem{notation}[theorem]{Notation}
\newtheorem{definition}[theorem]{Definition}

\theoremstyle{remark}
\newtheorem*{acknowledgments}{Acknowledgments}

\begin{document}

\title{Flabby strict deformation quantizations and $K$-groups}
\author{Hanfeng Li}

\address{Department of Mathematics \\
University of Toronto \\
Toronto ON M5S 3G3, CANADA} \email{hli@fields.toronto.edu}
\date{May 11, 2004}


\begin{abstract}
We construct examples of flabby strict deformation quantizations
not preserving $K$-groups. This answers a question of Rieffel
negatively.
\end{abstract}

\maketitle

\section{Introduction}
\label{intro:sec}

In the passage from classical mechanics to quantum mechanics, one
replaces smooth functions on symplectic manifolds (more generally,
Poisson manifolds) by operators on Hilbert spaces, and replaces
the Poisson bracket of smooth functions by commutators of
operators. Thinking of classical mechanics as limits of quantum
mechanics, one requires that the Poisson brackets becomes limits
of commutators.

There is an algebraic way of studying such process using formal
power series, called {\it deformation quantization} \cite{BFFLS78,
Weinstein95}. In order to study it in a stricter way, Rieffel
introduced \cite{Rieffel89} {\it strict deformation quantization}
of Poisson manifolds, within the framework of $C^*$-algebras. He
showed that noncommutative tori arise naturally as strict
deformation quantizations of the ordinary torus in the direction
of certain Poisson bracket. After that, a lot of interesting
examples of strict deformation quantizations have been
constructed. See \cite{Rieffel94, Rieffel98a} and the references
therein.

We refer the reader to \cite[Sections 10.1--10.3]{Dixmier77} for
the basic information about continuous fields of $C^*$-algebras.
Recall the definition of strict deformation quantization
\cite{Rieffel89, Rieffel98a}:
\begin{definition}\cite[Definition 1]{Rieffel98a} \label{sdq:def}
Let $M$ be a Poisson manifold, and let $C_{\infty}(M)$ be the
algebra of $\Ce$-valued continuous functions on $M$ vanishing at
$\infty$. By a \emph{strict deformation quantization} of M we mean
a dense $*$-subalgebra $A$ of $C_{\infty}(M)$ closed under the
Poisson bracket, together with a continuous field of
$C^*$-algebras $\mathcal{A}_{\hbar}$ over a closed subset $I$ of
the real line containing $0$ as a non-isolated point, and linear
maps $\pi_{\hbar}: A\to \mathcal{A}_{\hbar}$ for each $\hbar\in
I$, such that

(1)  $\mathcal{A}_0=C_{\infty}(M)$
 and $\pi_0$ is the canonical inclusion of
     $A$ into $C_{\infty}(M)$,

(2)  the section $(\pi_{\hbar}(f))$ is continuous for every $f\in
A$,

(3)  for all $f,g\in A$ we have
\begin{eqnarray*}
   \lim_{\hbar\to 0}\pa
   [\pi_{\hbar}(f), \, \pi_{\hbar}(g)]/(i\hbar)-\pi_{\hbar}(\{f,g\})\pa
   = 0,
\end{eqnarray*}

(4) $\pi_{\hbar}$ is injective and $\pi_{\hbar}(A)$ is a dense
$*$-subalgebra of $\mathcal{A}_{\hbar}$ for every $\hbar\in I$.

If $A\supseteq C^{\infty}_c(M)$, the space of compactly supported
$\Ce$-valued smooth functions on $M$, we say that the strict
deformation quantization is \emph{flabby}.
\end{definition}

Condition (4) above enables us to define a new $*$-algebra
structure and a new $C^*$-norm on $A$ at each $\hbar$ by pulling
back the $*$-algebra structure and norm of
$\pi_{\hbar}(A)\subseteq \mathcal{A}_{\hbar}$ to $A$ via
$\pi_{\hbar}$. Condition (2) means that this deformation of the
$*$-algebra structure and norm on $A$ is continuous.

Given a strict deformation quantization, a natural question is
whether the deformed $C^*$-algebras $\mathcal{A}_{\hbar}$ have the
same "algebraic topology", in particular, whether they have
isomorphic $K$-groups. Rieffel's quantization of Poisson manifolds
induced from actions of $\Re^d$ \cite{Rieffel93a} and many other
examples \cite{Nagy98} are known to preserve $K$-groups. Rieffel
showed examples of non-flabby strict deformation quantizations not
preserving $K$-groups, and asked \cite[Question 18]{Rieffel98a}:
{\it Are the $K$-groups of the deformed $C^*$-algebras of any
flabby strict deformation quantization all isomorphic?} A nice 
survey of various positive results on related problems may
be found in \cite{Rosenberg97}.

Shim \cite{Shim01} showed that above question has a negative answer
if one allows orbifolds. But it is not clear whether one can adapt
the method there to get smooth examples.

 Rieffel also pointed out \cite[page 321]{Rieffel98a} that in any strict
deformation quantization of a non-zero Poisson bracket if one
reparametrizes by replacing $\hbar$ by $\hbar^2$ one obtains a
strict deformation quantization of the $0$ Poisson bracket. Thus
to answer Rieffel's question it suffices to consider strict
deformation quantizations of the $0$ Poisson bracket.

The main purpose of this paper is to answer above question. In
Section~\ref{sdq:sec} we give a general method of constructing
flabby strict deformation quantization for the $0$ Poisson
bracket. In particular, we prove

\begin{theorem} \label{change K:theorem}
Let $M$ be a smooth manifold with $\dim M\ge 2$, equipped with the
$0$ Poisson bracket. If $\dim M$ is even (odd, resp.), then for any
integers $n_0\ge n_1 \ge 0$ ($n_1\ge n_0\ge 0$ resp.) there is a
flabby strict deformation quantization $\{\mathcal{A}_{\hbar},
\pi_{\hbar}\}_{\hbar\in I}$ of $M$ over $I=[0,1]$ with
$A=C^{\infty}_c(M)$ such that $K_i(\mathcal{A}_{\hbar})\cong
K_i(C_{\infty}(M))\oplus \Ze^{n_i}$ for all $0<\hbar\le 1$ and
$i=0,1$.
\end{theorem}

Theorem~\ref{change K:theorem} is far from being the most general
result one can obtain using our construction in
Section~\ref{sdq:sec}. However, it illustrates clearly that a lot
of manifolds equipped with the $0$ Poisson bracket have flabby
strict deformation quantizations not preserving $K$-groups.

In order to accommodate some other interesting examples such as
Berezin-Toeplitz quantization of K\"ahler manifolds,
Landsman introduced a weaker notion {\it strict quantization}
\cite[Definition II.1.1.1]{Landsman98b} \cite[Definition
23]{Rieffel98a}. This is defined in a way similar to a strict deformation
quantization, but without requiring the condition (4) in
Definition~\ref{sdq:def}. If $\pi_{\hbar}$ is injective for each
$\hbar\in I$ we say that the strict quantization is {\it
faithful}. It is natural to ask for the precise relation between
strict quantizations and strict deformation quantizations. Rieffel
also raised the question \cite[Question 25]{Rieffel98a}: {\it Is
there an example of a faithful strict quantization such that it is
impossible to restrict $\pi_{\hbar}$ to a dense $*$-subalgebra
$B\subseteq A$ to get a strict deformation quantization of $M$}?
Adapting our method in Section~\ref{sdq:sec} we also give such an
example for every manifold $M$ equipped with the $0$ Poisson
bracket. In \cite{Li3} strict quantizations are constructed for
every Poisson manifold, and it is impossible to restrict the
strict quantizations constructed there to dense $*$-subalgebras to
get strict deformation quantizations unless the Poisson bracket is
$0$ \cite[Corollary 5.6]{Li3}. Thus we get a complete answer to
Rieffel's question.

\begin{acknowledgments}
  I am grateful to Marc Rieffel for many helpful discussions and suggestions, 
and I also thank the referee for pointing out the reference \cite{Rosenberg97}.
\end{acknowledgments}

\section{Strict deformation quantizations for the $0$ Poisson bracket}
\label{sdq:sec}

We start with a general method of deforming a $C^*$-algebra. Let
$\mathcal{A}$ be a $C^*$-algebra and $A\subseteq \mathcal{A}$ a
dense $*$-subalgebra. Let $I(A)=\{b\in
\mathcal{M}(\mathcal{A}):bA, \, Ab\subseteq A\}$ be the idealizer
of $A$ in the multiplier algebra $\mathcal{M}(\mathcal{A})$ of
$\mathcal{A}$. Then $I(A)$ is a $*$-algebra containing $A$ as an
ideal, and for every $b\in (I(A))_{sa}$ clearly $bAb$ is a
$*$-subalgebra of $\mathcal{A}$. If furthermore the multiplication
by $b$ is injective on $A$, that is, $b\not \in Ann:=\{b'\in
\mathcal{M}(\mathcal{A}): b'a=0 \mbox{ for some } 0\neq a \in
A\}$, then we can pull back the multiplication and norm on $bAb$
to define a new multiplication $\times_b$ and a new norm $\pa
\cdot \pa_b$ on $A$ via the bijection $A\rightarrow bAb$.
Explicitly, $a\times_ba'=ab^2a'$ and $\pa a\pa_b=\pa bab\pa$. The
completion of $(A, \times_b, \pa \cdot \pa_b)$ is isomorphic to
$\overline{bAb}$ naturally.

Let $X$ be a topological space, and consider  a bounded map
$x\mapsto b_x$ from $X$ to $(I(A))_{sa}\subseteq
\mathcal{M}(\mathcal{A})$ continuous with respect to the strict
topology \cite[Definition 2.3.1]{Wegge93} on
$\mathcal{M}(\mathcal{A})$, i.e. the $\mathcal{A}$-valued
functions $x\mapsto \bar{a}\cdot b_x$ and $x\mapsto b_x \cdot
\bar{a}$ on $X$ are norm-continuous for every $\bar{a}\in
\mathcal{A}$. Then it follows easily that the $\mathcal{A}$-valued
function $x\mapsto b_xab_x$ on $X$ is norm-continuous for every
$a\in A$. Thus we get a continuous field of $C^*$-algebras
over $X$ with fibre algebra $\overline{b_xAb_x}$ at $x\in X$, as a
subfield of the trivial continuous field of $C^*$-algebras over
$X$ with fibres $\mathcal{A}$, and it contains $(b_xab_x)$ as a
continuous section for every $a\in A$.

Now we specialize to the commutative case. Let $M$ be a smooth
manifold, and let $\mathcal{A}=C_{\infty}(M), \,
A=C^{\infty}_c(M)$. Then $\mathcal{M}(\mathcal{A})$ is the space
$C_b(M)$ consisting of all $\Ce$-valued bounded continuous
functions on $M$, and the strict topology on $C_b(M)$ is
determined by uniform convergence on every compact subset of $M$.
The idealizer $I(A)$ is the space $C^{\infty}_b(M)$ consisting of
all $\Ce$-valued bounded smooth functions on $M$. Given $b \in
I(A)$, it is not in $Ann$ exactly if the zero set $Z_b$ of $b$ is
nowhere dense. Clearly $C_{\infty}(M\setminus Z_b)\supseteq bAb\supseteq C^{\infty}_c(M\setminus Z_b)$, and
hence $\overline{bAb}=C_{\infty}(M\setminus Z_b)$. Let
$X=I=[0,1]$. If $\hbar\mapsto b_{\hbar}$ is a bounded map from $I$
to $C^{\infty}_b(M)$ continuous with respect to the strict
topology on $C^{\infty}_b(M)$, then we get a continuous field of
$C^*$-algebras over $I$ with fibre $C_{\infty}(M\setminus Z_b)$ at
$\hbar$ and $(\pi_{\hbar}(a))$ is a continuous section for each
$a\in A$, where $\pi_{\hbar}(a)=b_{\hbar}ab_{\hbar}$. If
furthermore $b_0=1$ then the condition (1) of
Definition~\ref{sdq:def} is satisfied. Notice that when $M$ is
equipped with the $0$ Poisson bracket, the condition (3) of
Definition~\ref{sdq:def} holds trivially in our construction.
Summarizing above discussion we have reached:

\begin{proposition} \label{sdf non-unital:prop}
Let $M$ be a smooth manifold equipped with the $0$ Poisson
bracket. For any bounded map $\hbar\mapsto b_{\hbar}$ from
$I=[0,1]$ to $(C^{\infty}_b(M))_{sa}$ continuous respect to the
strict topology on $C^{\infty}_b(M)$, if $b_0=1$ and the zero set
$Z_{b_{\hbar}}$ of $b_{\hbar}$ is nowhere dense for every
$\hbar\in I$, then there is a flabby strict deformation
quantization $\{A_{\hbar}, \pi_{\hbar}\}_{\hbar\in I}$ of $M$ over
$I$ with $A=C^{\infty}_c(M)$ and
$\mathcal{A}_{\hbar}=C_{\infty}(M\setminus Z_{b_{\hbar}})$ for
every $\hbar \in I$.
\end{proposition}

\begin{example} \label{R^n:example}
Let $M=\Re^n$. Take a bounded smooth real-valued function $F$ on
$\Re^n$ such that $F=1$ in a neighborhood of the origin and $F$
vanishes exactly at one point $P$. Set $b_{\hbar}(x)=F(\hbar x)$
for all $0\le \hbar\le 1$ and $x\in \Re^n$. Then for each
$0<\hbar\le 1$ the space $\Re^n\setminus
Z_{b_{\hbar}}=\Re^n\setminus \{P/\hbar\}$ is homeomorphic to
$\Re\times S^{n-1}$. Now by Proposition~\ref{sdf non-unital:prop}
there is a flabby strict deformation quantization $\{A_{\hbar},
\pi_{\hbar}\}_{\hbar \in I}$ of $\Re^n$ over $I=[0,1]$ with
$A=C^{\infty}_c(\Re^n)$ and
$\mathcal{A}_{\hbar}=C_{\infty}(M\setminus Z_{b_{\hbar}})\cong
C_{\infty}(\Re \times S^{n-1})$ for every $0< \hbar \le 1$. Then
by the Bott periodicity $K_i(\mathcal{A}_{\hbar})\cong
K_{i+1}(C(S^{n-1}))$ \cite[Theorem 7.2.5, page 158]{Wegge93} for
every $0<\hbar\le 1$ and $i=0, 1$. Thus when $n$ is odd
$K_i(\mathcal{A}_{\hbar})\not \cong K_i(C_{\infty}(\Re^n))$ for
all $0<\hbar\le 1$ and $i=0,1$ (see for instance \cite[page
123]{Wegge93} for the $K$-groups of $\Re^n$ and $S^n$). When $n$
is even, $K_1(\mathcal{A}_{\hbar})\not \cong
K_1(C_{\infty}(\Re^n))$ for all $0<\hbar\le 1$.
\end{example}

When $M$ is compact, in Proposition~\ref{sdf non-unital:prop} the
element $b_{\hbar}$ has to be invertible in $C(M)$ for small
$\hbar$ and consequently $\mathcal{A}_{\hbar}=C(M)$. Thus in order
to construct strict deformation quantizations for compact $M$ such
that the $K$-groups of $\mathcal{A}_{\hbar}$ are not isomorphic to
those of $C(M)$ for any $\hbar\neq 0$, we have to modify the
construction in Proposition~\ref{sdf non-unital:prop}. Notice that
if we set $\pi'_{\hbar}(a+\lambda)=b_{\hbar}ab_{\hbar}+\lambda$
for $a\in C^{\infty}_c(\Re^n), \, \lambda \in \Ce$ in
Example~\ref{R^n:example}, then we get a strict deformation
quantization of $S^n$ equipped with the $0$ Poisson bracket. This
leads to Proposition~\ref{sdf unital:prop} below.

\begin{notation} \label{Fm:notation}
We denote by $\mathcal{F}_m$ the space of smooth real-valued
functions $F$ on $\Re^m$ such that $F$ is equal to $1$ outside a
compact subset of $\Re^m$ and the zero set $Z_F$ of $F$ is nowhere
dense.
\end{notation}

\begin{proposition}  \label{sdf unital:prop}
Let $M$ be a smooth manifold equipped with the $0$ Poisson
bracket. Let $U$ be an open subset of $M$ with a diffeomorphism
$\varphi:U\rightarrow \Re^m$. For any $F\in \mathcal{F}_m$ there
is a flabby strict deformation quantization $\{A_{\hbar},
\pi_{\hbar}\}_{\hbar \in I}$ of $M$ over $I=[0,1]$ with
$A=C^{\infty}_c(M)$ such that $\mathcal{A}_{\hbar}\cong
C_{\infty}(M/ Y)$ for every $0<\hbar \le 1$, where
$Y=\varphi^{-1}(Z_F\cup \{0\})$.
\end{proposition}
\begin{proof}
Set $F_0=1$ and $F_{\hbar}(x)=F(x/\hbar)$ for all $0<\hbar\le 1$
and $x\in \Re^m$. Then $F_{\hbar}\in \mathcal{F}_m$ for each
$\hbar\in I$ and we can extend the pull-back $F_{\hbar}\circ
\varphi \in C^{\infty}(U)$ to a smooth function $b_{\hbar}$ on $M$
by setting it to be $1$ outside $U$. Clearly
$b_{\hbar}A'b_{\hbar}$ is a $*$-subalgebra of $A'$. Notice that
there is a compact set $W\subset U$ such that $b_{\hbar}=1$ on
$M\setminus W$ for all $\hbar\in I$, and $W$ contains $\varphi^{-1}(0)$. 
Take an $H\in
(C^{\infty}_c(M))_{\Re}$ such that $H=1$ on $W$. 
Denote by $A'$ the space of functions in $C^{\infty}_c(M)$
vanishing at $\varphi^{-1}(0)$. Then $C^{\infty}_c(M)=A'\oplus \Ce
H$ as complex vector spaces, and 
$H^2-H=b_{\hbar}(H^2-H)b_{\hbar}\in b_{\hbar}A'b_{\hbar}$. It is
easy to see that $b_{\hbar}A'b_{\hbar}+\Ce H$ is a $*$-subalgebra
of $C^{\infty}_c(M)$ and the linear map
$\pi_{\hbar}:C^{\infty}_c(M)\rightarrow b_{\hbar}A'b_{\hbar}+\Ce
H$ defined by $\pi_{\hbar}(a'+\lambda
H)=b_{\hbar}a'b_{\hbar}+\lambda H$ for $a'\in A'$ and $\lambda\in
\Ce$ is bijective. For each $a'\in A'$ clearly the map
$\hbar\mapsto b_{\hbar}a'b_{\hbar}\in C_{\infty}(M)$ is continuous
on $I=[0,1]$. Thus for each $a\in A=C^{\infty}_c(M)$,
$(\pi_{\hbar}(a))$ is a continuous section in the continuous
subfield $\{\mathcal{A}_{\hbar}=\overline{b_{\hbar}A'b_{\hbar}+\Ce
H}\}_{\hbar\in I}$ of the trivial field of $C^*$-algebras over $I$
with fibres $C_{\infty}(M)$. Therefore $\{\mathcal{A}_{\hbar},
\pi_{\hbar}\}_{\hbar \in I}$ is a flabby strict deformation
quantization of $M$.

Set $Y_{\hbar}=\varphi^{-1}((\hbar Z_F)\cup \{0\})$. Clearly
$C_{\infty}(M\setminus Y_{\hbar})\supseteq b_{\hbar}A'b_{\hbar}\supseteq C^{\infty}_c(M\setminus Y_{\hbar})$. Thus
$\overline{b_{\hbar}A'b_{\hbar}+\Ce
H}=\overline{b_{\hbar}A'b_{\hbar}}+\Ce H$ is exactly the space of
functions in $C_{\infty}(M)$ taking the same value on $Y_{\hbar}$,
which is just $C_{\infty}(M/ Y_{\hbar})$. When $0<\hbar \le 1$,
the space $M/ Y_{\hbar}$ is homeomorphic to $M/Y$, and hence
$\mathcal{A}_{\hbar}=\overline{b_{\hbar}A'b_{\hbar}+\Ce H}\cong
C_{\infty}(M/Y)$ as desired.
\end{proof}

Next we describe a case in which we can relate the $K$-groups of
$C_{\infty}(M/Y)$ to those of $C_{\infty}(M)$ easily:

\begin{lemma} \label{K:lemma}
Let $D$ be the subset of $\Re^m$ consisting of points $(x_1,
{\cdots}, x_m)$ with $0< x_1, {\cdots}, x_m< 1$. Let $M, \varphi,
F$ and $Y$ be as in Proposition~\ref{sdf unital:prop}. Suppose
that $\partial D\subseteq Z_F\subseteq \bar{D}$. Then
\begin{eqnarray*}
K_i(C_{\infty}(M\setminus Y))\cong K_i(C_{\infty}(M))\oplus
K_i(C_{\infty}(D\setminus Z_F))
\end{eqnarray*}
for $i=0,1$.
\end{lemma}
\begin{proof}
Let $\phi:M\rightarrow M/Y$ be the quotient map, and let
$W=\phi(M\setminus \varphi^{-1}(D))$. Then $W$ is a closed subset
of $M/Y$, and the complement is homeomorphic to $D\setminus Z_F$.
Define a map $\psi:M/Y\rightarrow W$ as the identity map on $W$
and $\psi((M/F)\setminus W)=\phi(Y)$. Then $\psi$ is continuous
and proper, i.e. the inverse image of every compact subset of $W$
is compact. Thus the exact sequence
\begin{eqnarray*}
0\rightarrow C_{\infty}(D\setminus Z_F)\rightarrow
C_{\infty}(M/Y)\rightarrow C_{\infty}(W)\rightarrow 0
\end{eqnarray*}
splits. Therefore $K_i(C_{\infty}(M\setminus Y))\cong
K_i(C_{\infty}(W))\oplus K_i(C_{\infty}(D\setminus Z_F))$ for
$i=0,1$. Now Lemma~\ref{K:lemma} follows from the fact that $W$ is
homeomorphic to $M$.
\end{proof}

Notice that if a compact set $Z\subseteq \Re^m$ is the zero set of
some non-negative $f\in C^{\infty}(M)$, then it is also the zero
set of some $F\in \mathcal{F}_m$ (for instance, take a
non-negative $g\in C^{\infty}_c(M)$ with $g|_Z=1$ and set
$F(x)=\frac{f(x)}{f(x)+g(x)}$ for all $x\in \Re^m$). Also notice
that if closed subsets $Z_1$ and $Z_2$ of $\Re^m$ are the zero
sets of non-negative smooth functions on $\Re^m$, then so are
$Z_1\cap Z_2$ and $Z_1\cup Z_2$. From these observation we get
easily
\begin{lemma} \label{function:lemma}
Let $m\ge 2$, and let $D$ be as in Lemma~\ref{K:lemma}. For any
$k, j\ge 0$ and $1\le s \le 2$ there exits an $F\in \mathcal{F}_m$
satisfying  $\partial D\subseteq Z_F\subseteq \bar{D}$ such that
$D\setminus Z_F$ is homeomorphic to the disjoint union of $k$ many
$\Re^m$ and $j$ many $\Re^s\times S^{m-s}$.
\end{lemma}

Now we are ready to prove Theorem~\ref{change K:theorem}:
\begin{proof}[Proof of Theorem~\ref{change K:theorem}]
The case in which $\dim M$ is even follows from
Proposition~\ref{sdf unital:prop} and Lemma~\ref{K:lemma} by
taking $k=n_0-n_1, \, j=n_1,\, s=1$ in Lemma~\ref{function:lemma}.
Similarly the case in which $\dim M$ is odd follows by taking
$k=n_1-n_0,\, j=n_0,\, s=2$.
\end{proof}

Finally we discuss how to adapt our method to construct strict
quantizations which can't be restricted to dense $*$-subalgebras
to yield strict deformation quantizations. Notice that if we relax
the condition $b_{\hbar}\in (C^{\infty}_b(M))_{sa}$ in
Proposition~\ref{sdf non-unital:prop} to $b_{\hbar}\in
(C_b(M))_{sa}$ and set $\mathcal{A}_{\hbar}$ to be the
$C^*$-subalgebra of $_{\infty}(M)$ generated by
$b_{\hbar}Ab_{\hbar}$, then we get a faithful strict quantization
of $M$ over $I=[0,1]$ with $A=C^{\infty}_c(M)$ and
$\pi_{\hbar}(a)=b_{\hbar}ab_{\hbar}$. Take a nonnegative $F\in
C_b(M)$ such that $F$ is not smooth at some point $P$. Set
$b_{\hbar}=((1-\hbar)+\hbar F)^{1/2}$ for every $\hbar \in I$. Let
$B\subseteq A$ be a dense $*$-subalgebra. Then we can find $f\in
B$ such that $f(P)\neq 0$. Clearly $f^2b^2_{\hbar}$ is not smooth
at $P$ for $0<\hbar\le 1$. It follows that $(\pi_{\hbar}(f))^2$ is
not in $\pi_{\hbar}(B)$ for $0<\hbar\le 1$. Thus we get:
\begin{proposition} \label{sq:prop}
Let $M$ be a smooth manifold equipped with the $0$ Poisson
bracket. Then there is a faithful strict quantization
$\{\mathcal{A}_{\hbar}, \pi_{\hbar}\}_{\hbar\in I}$ of $M$ over
$I=[0,1]$ with $A=C^{\infty}_c(M)$ such that it is impossible to
restrict $\pi_{\hbar}$ to a dense $*$-subalgebra $B\subseteq A$ to
get a strict deformation quantization of $M$.
\end{proposition}


\end{document}